\documentclass[11pt, reqno]{amsart}
\usepackage{graphics}
\usepackage{amsfonts}

\newcommand{\bproof}{{\bf\underline{Proof:}} }
\def\Box{\vcenter{\vbox{\hrule\hbox{\vrule
     \vbox to 8.8pt{\hbox to 10pt{}\vfill}\vrule}\hrule}}}

\newcommand{\Ff}{{\mathbb F}}
\newcommand{\Zz}{{\mathbb Z}}

\newcommand{\eproof}{\hfill$\Box$\vspace{4mm}}

\newcommand{\Z}{{\mathbb Z}}

\newtheorem{thm}{Theorem}[section]
\newtheorem{lemma}[thm]{Lemma}
\newtheorem{cor}[thm]{Corollary}

\newtheorem{example}[thm]{Example}
\newtheorem{remark}[thm]{Remark}
\def\Tr{\operatorname{Tr}}
\def\tr{\operatorname{tr}}
\DeclareMathOperator{\PG}{PG} \DeclareMathOperator{\GR}{GR}

\numberwithin{equation}{section}

\begin{document}

\title[amorphic association schemes]
{Amorphic Association Schemes with Negative Latin Square Type Graphs}

\author{James A. Davis, Qing Xiang}

\dedicatory{Dedicated to Zhe-Xian Wan on the occasion of his 80th birthday}

\address{Department of Mathematics and Computer Science, University of Richmond, Richmond, VA 23173, USA,
email: {\tt jdavis@richmond.edu}}

\address{Department of Mathematical Sciences, University of Delaware, Newark, DE 19716, USA,
email: {\tt xiang@math.udel.edu}}


\date{}

\begin{abstract}
Applying results from partial difference sets, quadratic forms, and recent
results of Brouwer and Van Dam, we construct the first
known amorphic association scheme with negative Latin square type
graphs and whose underlying set is a nonelementary abelian 2-group.
We give a simple proof of a result of Hamilton that generalizes
Brouwer's result.  We use multiple distinct quadratic forms to construct
amorphic association schemes with a large number of classes.
\end{abstract}
\maketitle

\section{Introduction}
\label{intro}

Let $X$ be a finite set. A ({\it symmetric}) {\it association scheme} with $d$ classes on $X$ is a partition of
$X\times X$ into sets $R_0$, $R_1, \ldots , R_d$ (called relations, or associate classes) such that

\begin{enumerate}
\item $R_0=\{(x,x)\mid x\in X\}$ (the diagonal relation);
\item $R_{\ell}$ is symmetric for $\ell=1,2,\ldots ,d$;
\item for all $i,j,k$ in $\{0,1,2,\ldots ,d\}$ there is an integer $p_{ij}^k$ such that, for all $(x,y)\in R_k$,
$$|\{z\in X \mid (x,z)\in R_i\; {\rm and}\; (z,y)\in R_j\}|=p_{ij}^k.$$
\end{enumerate}

Since each symmetric relation $R_\ell$, $1\leq \ell\leq d$, corresponds to an (undirected) graph $G_\ell=(X,
R_{\ell})$, $1\leq \ell\leq d$, with vertex set $X$ and edge set $R_{\ell}$, we can think of an association
scheme $(X, \{R_\ell\}_{0\leq\ell\leq d})$ as an edge-decomposition of the complete graph on the vertex set $X$
into graphs $G_\ell$ on the same vertex set with the property that for all $i,j,k$ in $\{1,2,\ldots ,d\}$ and for
all $xy\in E(G_k)$,
$$|\{z\in X \mid xz\in E(G_i)\; {\rm and}\; zy\in E(G_j)\}|=p_{ij}^k,$$
where $E(G_k)$, $E(G_i)$, and $E(G_j)$ are the edge sets of $G_k$, $G_i$ and $G_j$ respectively. The graphs
$G_{\ell}$, $1\leq\ell\leq d$, will be called {\it the graphs} of the association scheme $(X, \{R_\ell\}_{0\leq \ell\leq d})$. (Note
that it follows from the definition of an association scheme that each graph $G_{\ell}$ of the association scheme
$(X, \{R_\ell\}_{0\leq \ell\leq d})$ is regular with valency $n_{\ell}=p_{\ell \ell}^0$.) A {\it strongly regular
graph} (SRG) is an association scheme with two classes, say $(X,\{R_{\ell}\}_{0\leq \ell \leq 2})$. The elements
of $X$ are the vertices of the graph and $\{x,y\}$ is an edge if $(x,y)\in R_1$. There are four parameters
$(v,k,\lambda,\mu)$ associated with the SRG, where
$$v=|X|, \; k=n_1,\; \lambda=p_{11}^1,\; {\rm and}\; \mu=p_{11}^2.$$
For more background on association schemes and strongly regular graphs, see~\cite{bcn} and~\cite{vanLintWilson}.
Given an association scheme $(X, \{R_\ell\}_{0\leq \ell\leq d})$, we can take the union of classes to form graphs
with larger edge sets (this is called a {\em{fusion}}), but it is not necessarily guaranteed that the fused
collection of graphs will form an association scheme on $X$. If an association scheme has the property that any
of its fusions is also an association scheme, then we call the association scheme {\em{amorphic}}. A typical
example of amorphic association schemes is given by the so-called uniform cyclotomy of finite fields.

\begin{example}\label{cycexamp}
Let $q=p^s$ be a power of a prime $p$ and let $q-1=ef$ with $e>1$. Let $\Ff_q^*$ denote the multiplicative group
of $\Ff_q$, $C_0$ be the subgroup of $\Ff_q^*$ of index $e$, and let $C_0,C_1,\ldots ,C_{e-1}$ be the cosets of
$C_0$ in $\Ff_q^*$. Assume $-1\in C_0$. Define $R_0=\{(x,x)\mid  x\in \Ff_q\}$, and for $i\in \{1,2,\ldots ,e\}$,
define $R_i=\{(x,y)\mid x,y\in \Ff_q, x-y\in C_{i-1}\}$. Then $(\Ff_q, \{R_i\}_{0\leq i\leq e})$ is an
association scheme with $e$ classes. It was proved in \cite{bmw} that, for $e>2$, this scheme is amorphic if and
only if $-1$ is a power of $p$ modulo $e$.
\end{example}

A $(v,k,\lambda,\mu)$ strongly regular graph is said to be of {\em{Latin square type}} (resp. {\it negative Latin
square type}) if $(v,k,\lambda,\mu) = (n^2, r(n-\epsilon), \epsilon n+r^2-3 \epsilon r, r^2 - \epsilon r)$ and
$\epsilon = 1$ (resp. $\epsilon=-1$). If an association scheme is amorphic, then each of its graphs is clearly
strongly regular. Moreover, A. V. Ivanov \cite{ivanov} showed that in an amorphic association scheme with at
least three classes all graphs of the scheme are of Latin square type, or all graphs are of negative Latin square
type. The converse of Ivanov's result is proved to be true in \cite{imy}. In fact even more is true because Van
Dam~\cite{vanDam} could prove the following result.

\begin{thm}
\label{vanDamtheorem} Let $X$ be a set of size $v$, let $\{G_1, G_2, \ldots, G_d \}$ be an edge-decomposition of
the complete graph on $X$, where each $G_i$ is a strongly regular graph on $X$. If $G_i, 1 \leq i \leq d,$ are
all of Latin square type or all of negative Latin square type, then the decomposition is a $d$-class amorphic
association scheme on $X$.
\end{thm}

We will consider the following situation. Let $G$ be a finite additive group with identity $0$. If we can
partition $G\setminus \{0\}$ into sets $D_i$, all of which are Latin square type partial difference sets (defined
below), or all of which are negative Latin square type partial difference sets, then the corresponding strongly
regular Cayley graphs will satisfy the conditions of Theorem~\ref{vanDamtheorem} and hence will form an amorphic
association scheme. That is our objective in the constructions to follow.

\subsection{Partial Difference Sets}
A $k$-element subset $D$ of a finite multiplicative group $G$ of order $v$ is called a $(v,k,\lambda,
\mu)$-{\em{partial difference set}} (PDS) in $G$ provided that the list of ``differences'',  $d_1 d_2^{-1}$,
$d_1, d_2\in D$, $d_1 \neq d_2,$ contains each nonidentity element of $D$ exactly $\lambda$ times and each
nonidentity element in $G \backslash D$ exactly $\mu$ times. Partial difference sets are equivalent to strongly
regular graphs with a regular automorphism group via the Cayley graph construction. We refer the reader to
Ma~\cite{masurvey} for a survey of results on PDS. If a PDS has parameters $(n^2, r(n-\epsilon), \epsilon n+r^2-3
\epsilon r, r^2 - \epsilon r)$ with $\epsilon = 1$ (resp. $\epsilon=-1$), then the PDS is said to be of {\it
Latin square type} (resp. {\it negative Latin square type}). Latin square type PDSs have been constructed in a
variety of nonisomorphic groups of order $n^2$ by taking the disjoint union of $r$ subgroups of order $n$ that
pairwise intersect trivially (\cite{masurvey}). In contrast, most known constructions of negative Latin square
type PDSs occur in elementary abelian $p$-groups. As far as we know, the first infinite family of negative Latin
square type PDSs in nonelementary abelian $p$-groups was constructed recently in \cite{davisxiang}.

A complex character of an abelian group is a homomorphism from the group to the multiplicative group of complex
roots of unity. The {\em principal character} is the character mapping every element of the group to 1.  All
other characters are called {\it nonprincipal}. Starting with the important work of Turyn \cite{turyn}, character
sums have been a powerful tool in the study of difference sets of all types.  The following lemma states how
character sums can be used to verify that a subset of an abelian group is a PDS.

\begin{lemma}\label{latinPDS}
\label{charsum-PDS} Let $G$ be an abelian group of order $v$ and
$D$ be a $k$-subset of $G$ such that $\{ d^{-1} \mid d \in D \} =
D$ and $1\not\in D$. Then $D$ is a $(v,k,\lambda, \mu)$-PDS in $G$
if and only if, for any complex character $\chi$ of $G$,
$$\sum_{d \in D} \chi(d) = \left\{ \begin{array}{ll}
    k, & \mbox{if $\chi$ is principal on $G$,} \\
    \frac{(\lambda - \mu) \pm \sqrt{(\lambda - \mu)^2 + 4(k-\mu)}}{2}, & \mbox{if $\chi$ is nonprincipal on $G$.}
            \end{array}
        \right. $$
In particular, let $v=n^2$ and let $k=r(n-\epsilon)$, where $\epsilon=\pm 1$. Then $D$ is a $(n^2, r(n-\epsilon),
\epsilon n+r^2-3\epsilon r, r^2-\epsilon r)$-PDS in $G$ if and only if, for any complex character $\chi$ of $G$,
$$\sum_{d \in D} \chi(d) = \left\{ \begin{array}{ll}
    r(n-\epsilon), & \mbox{if $\chi$ is principal on $G$,} \\
    \epsilon (n-r)\; {\mbox or}\; -\epsilon r, & \mbox{if $\chi$ is nonprincipal on $G$.}
            \end{array}
        \right. $$
\end{lemma}

We will use Lemma~\ref{latinPDS} frequently in the rest of the paper.

\subsection{Galois Ring Preliminaries}

Interested readers are referred to \cite{mcdonald} for more details on Galois rings.  We will follow the approach
outlined in~\cite{davisxiang}, including the fact that we will only use Galois rings over $\Z_4$. A {\em Galois
ring over $\Z_{4}$ of degree $t$}, $t\geq 2$, denoted $\GR(4,t)$, is the quotient ring $\Z_{4}[x]/ \langle
\Phi(x) \rangle$, where $\Phi(x)$ is a basic primitive polynomial in $\Z_{4}[x]$ of degree $t$. Let $\xi$ be a
root of $\Phi(x)$ in $\GR(4,t)$. Then we have $\GR(4,t)= \Z_4[\xi]$.  In this paper, we will only need $\GR(4,2) =
\Z_4[x]/ \langle x^2+x+1 \rangle$.

The ring $R=\GR(4,2)$ is a finite local ring with unique maximal
ideal $2R$, and $R/2R$ is isomorphic to the finite field
$\Ff_{4}$. If we denote the natural epimorphism from $R$ to
$R/2R\cong \Ff_{4}$ by $\pi$, then $\alpha = \pi(\xi)$ is a
primitive element of $\Ff_{4}$.

The set ${\mathcal T} = \{ 0,1,\xi,\xi^2 \}$ is a complete set of
coset representatives of $2R$ in $R$. This set is usually called a
{\em Teichm\"uller system} of $R$.  The restriction of $\pi$ to
${\mathcal T}$ is a bijection from ${\mathcal T}$ to $\Ff_{4}$,
and we refer to this bijection as $\pi_{_{\mathcal T}}$. An
arbitrary element $\beta$ of $R$ has a unique 2-adic
representation
$$\beta = \beta_1 + 2 \beta_2,$$
where $\beta_1, \beta_2 \in {\mathcal T}$.  Combining $\pi_{_{\mathcal T}}$ with this 2-adic representation, we
get a bijection $F$ from $\Ff_4^{2 \ell}$ to $R\times \Ff_4^{2 \ell - 2}$ defined by
\begin{eqnarray*}
F: (x_1, x_2, \ldots,  x_{2 \ell}) \mapsto & (\pi_{_{\mathcal
T}}^{-1}(x_1) + 2 \pi_{_{\mathcal T}}^{-1}(x_2), x_3, x_4, \ldots,
 x_{2 \ell}).\\
\end{eqnarray*}
The inverse of this map is the map $F^{-1}$ from $R \times \Ff_4^{2 \ell - 2}$ to $\Ff_4^{2 \ell}$,
\begin{eqnarray*}
F^{-1}: (\xi_1+2\xi_2,x_3, \ldots , x_{2\ell})\mapsto
&(\pi_{_{\mathcal T}}(\xi_1),\pi_{_{\mathcal T}}(\xi_2), x_3,
\ldots
, x_{2\ell}).\\
\end{eqnarray*}
To simplify notation we will usually omit the subindex in the bijection $\pi_{_{\mathcal
T}}^{-1}:\Ff_{4}\rightarrow {\mathcal T}$. (So from now on, $\pi^{-1}$ means the inverse of the bijection
$\pi_{_{\mathcal T}}:{\mathcal T}\rightarrow \Ff_{4}$.) We will show in
Section~\ref{nonelementaryabelianconstruction} how we can use $F$ as a character-sum-preserving bijection to
construct a PDS in a nonelementary abelian group, namely the additive group of $R\times \Ff_4^{2 \ell - 2}$.

The {\em Frobenius map} $f$ from $R$ to itself is the ring automorphism $f: \beta_1+2\beta_2 \mapsto \beta_1^2 +
2 \beta_2^2$. This map is used to define the {\em trace} ${\rm Tr}$ from $R$ to $\Z_{4}$, namely, ${\rm
Tr}(\beta) = \beta + \beta^f$, for $\beta\in R$.  We note that here the Galois ring trace $\Tr: R\rightarrow
\Zz_4$ is related to the finite field trace $\tr:\Ff_{4}\rightarrow \Ff_2$ via
\begin{equation}\label{commutrace}
\tr\circ \pi=\pi\circ \Tr. \end{equation} As a consequence, we have
$$\sqrt{-1}^{\Tr(2x)} =\sqrt{-1}^{2\Tr(x)}=(-1)^{\Tr(x)}=(-1)^{\pi\circ \Tr(x)}=(-1)^{\tr(\pi(x))},$$
for all $x\in R$.

All additive characters of Galois rings and finite fields can be defined by using the appropriate trace, as
indicated in the following well-known lemma.

\begin{lemma}
\label{tracegiveschar} \begin{enumerate}
\item  Let $\chi$ be
an additive character of $\Ff_q$, where $q$ is a power of $p$, $p$ is a prime,
and let $\eta_p$ be a complex primitive $p^{\rm th}$ root of unity. Then there
is a $w \in \Ff_q$ such that $\chi(x) = \eta_p^{{\rm tr}(w x)}$ for all $x \in \Ff_q$, where ${\rm tr}$ is the trace from $\Ff_q$ to $\Ff_p$.

\item Let $\psi$ be an additive character of $\GR(4,t)$. Then there is a
$\beta \in \GR(4,t)$ such that $\psi(x) = \sqrt{-1}^{{\rm Tr}(\beta x)}$ for all $x \in \GR(4,t)$.
\end{enumerate}
\end{lemma}

Since we can write $\beta = \beta_1 + 2 \beta_2$ for any $\beta \in R=\GR(4,2)$, where $\beta_k \in {\mathcal T},
k=1,2$, we will use the notation $\psi_{_{\beta}} = \psi_{_{\beta_1 + 2 \beta_2}}$ indicating the ring element
used to define the character $x\mapsto \sqrt{-1}^{{\rm Tr}(\beta x)}$. If $\beta_1 =0$ but $\beta_2\neq 0$, then
$\psi_{_{2 \beta_2}}$ is a character of order 2 and $\psi_{_{2 \beta_2}}$ is principal on $2R$. If $\beta_1 \neq
0$, then $\psi_{_{\beta_1 + 2 \beta_2}}$ is a character of order 4 and $\psi_{_{\beta_1 + 2 \beta_2}}$ is
nonprincipal on $2R$.  For the finite field $\Ff_{4}$, we will use the notation $\chi_w$ to indicate the field
element used to define the character $x\mapsto (-1)^{{\rm tr}(w x)}$. Characters of $R\times \Ff_4^{2 \ell - 2}$
will be written as
\begin{equation}\label{defchar}
\Psi = \psi_{_{\beta_{1} + 2 \beta_{2}}} \otimes \chi_{_{w_{3}}} \otimes \chi_{_{w_{4}}} \otimes \cdots \otimes
\chi_{_{w_{2 \ell}}} = \psi_{_{\beta_{1} + 2 \beta_{2}}} \otimes \chi_{_{(w_{3}, w_4, \ldots, w_{2 \ell})}},
\end{equation}
where $\beta_i \in {\mathcal T}, 1 \leq i \leq 2$, and $w_{i} \in \Ff_4, 3 \leq i \leq 2 \ell$.

\subsection{Quadratic Forms}
\label{quadraticforms}

Let $\Ff_q$ be the field of order $q$, where $q$ is a prime power, and let $V$ be an $m$-dimensional vector space
over $\Ff_q$. A function $Q: V \rightarrow \Ff_q$ is called a {\it{quadratic form}} if
\begin{enumerate}
\item  $Q(\gamma v) = \gamma^2 Q(v)$ for all $\gamma \in
\Ff_q$ and
$v \in V$,
\item the function $B: V \times V \rightarrow \Ff_q$ defined by
$B(v_1,v_2) = Q(v_1+v_2) - Q(v_1) - Q(v_2)$ is bilinear.
\end{enumerate}
We call $Q$ {\em{nonsingular}} if the subspace $W$ with the property that $Q$ vanishes on $W$ and $B(w,v)=0$ for
all $v\in V$ and $w\in W$ is the zero subspace. For background on quadratic forms, see~\cite{cameron}.  Let $Q$
be a nonsingular quadratic form on an $m$-dimensional vector space $V$ over $\Ff_q$. If $m$ is even, then $Q$ is
equivalent to either $x_1x_2 + x_3 x_4 + \cdots + x_{m-1} x_{m}$ (called {\em{hyperbolic}}, or type $+1$) or
$p(x_1, x_2) + x_3 x_4 + \cdots + x_{m-3}x_{m-2} + x_{m-1}x_{m}$, where $p(x_1,x_2)$ is an irreducible quadratic
form in two indeterminates (called {\em{elliptic}}, or type $-1$).

Let $\PG(m-1,q)$ be the $(m-1)$-dimensional desarguesian projective space over $\Ff_q$, and let $V$ be its
associated vector space. The {\it quadric} of the projective space $\PG(m-1,q)$ defined by a quadratic form
$Q: V\rightarrow \Ff_q$ is the point set ${\mathcal Q}=\{\langle v\rangle\in \PG(m-1,q)\mid Q(v)=0\}$. The
following theorem about quadrics in $\PG(m-1,q)$ is well known, see for example \cite{calderbankkantor}.

\begin{thm}  \label{ellquadform}
Let ${\mathcal Q}$ be a nonsingular elliptic quadric in $\PG(2\ell-1,q)$, and let $\Omega$ denote the set of
nonzero vectors in $\Ff_q^{2\ell}$ corresponding to ${\mathcal Q}$ (i.e., $\Omega=\{x\in \Ff_q^{2\ell}\mid x\neq
0, \langle x\rangle\in {\mathcal Q}\}$). For any nontrivial additive character $\chi$ of  $\Ff_q^{2\ell}$, we have
$$\chi(\Omega) = \left\{ \begin{array}{ll}
    (q^{\ell-1}-1)-q^{\ell}, & \mbox{or} \\
    (q^{\ell-1}-1).
            \end{array}
        \right. $$
That is, $\Omega$ is a
$(q^{2\ell},(q^{\ell}+1)(q^{\ell-1}-1),q^{2\ell-2}-q^{\ell-1}(q-1)-2,q^{2\ell-2}-q^{\ell-1})$-negative Latin
square type PDS in the additive group of $\Ff_q^{2\ell}$.
\end{thm}

There is a corresponding theorem for hyperbolic quadratic forms,
where the set $\Omega$ in that case will be a Latin square type
PDS.  Since the main focus of this paper is on constructions of
negative Latin square type PDSs, we do not state that theorem nor
do we include the constructions of Latin square type amorphic
association schemes.

Brouwer~\cite{brouwer} showed that one can also use quadratic forms to define other PDSs than the ones listed in
the previous theorem. Suppose $\Ff_{q_{_{0}}}$ is a subfield of $\Ff_q$, say $q=q_{_{0}}^e$. Let
$V=\Ff_q^{2\ell}$, $Q: V\rightarrow \Ff_q$ be a nonsingular quadratic form on $V$, and $\tr_{q/q_{_{0}}}$ be the
trace from $\Ff_q$ to $\Ff_{q_{_{0}}}$. Then $Q_0=\tr_{q/q_{_{0}}}\circ Q: V\rightarrow \Ff_{q_{_{0}}}$ is a
quadratic form on the same $V$ but now viewed as a $2\ell e$-dimensional vector space over $\Ff_{q_{_{0}}}$.
Define
$$\Omega=\{x\in V\mid Q(x)=0, x\neq 0\}, \;
\Omega_0=\{x\in V\mid Q_0(x)=0, x\neq 0\}.$$ Then clearly $\Omega\subset \Omega_0$. Brouwer~\cite{brouwer} proved
the following theorem.

\begin{thm}\label{brouwerthm}
With the assumptions and notation above, the set $\Omega_0\setminus\Omega$ is a PDS in the additive group of $V$.
\end{thm}

\begin{remark}
In the above theorem, $Q$ can be either elliptic or hyperbolic. If $Q$ is a nonsingular elliptic quadratic form
on $V=\Ff_q^{2\ell}$, then
$$|\Omega_0\setminus\Omega|=(q_{_{0}}^{\ell e}+1)(q_{_{0}}^{\ell e-1}-q_{_{0}}^{\ell e-e}),$$
and $\Omega_0\setminus\Omega$ is a negative Latin square type PDS in $(V, +)$. This will be the case to which we
pay most attention. If $Q$ is a nonsingular hyperbolic quadratic form on $V=\Ff_q^{2\ell}$, then
$$|\Omega_0\setminus\Omega|=(q_{_{0}}^{\ell e}-1)(q_{_{0}}^{\ell e-1}-q_{_{0}}^{\ell e-e}),$$
and $\Omega_0\setminus\Omega$ is a Latin square type PDS in $(V, +)$.

\end{remark}

We will see in the following sections how this result leads to amorphic association schemes.  In
Section~\ref{example}, we provide a relatively simple example that motivated the work done in this paper.
Section~\ref{nonelementaryabelianconstruction} uses techniques similar to that in \cite{davisxiang} to construct
an amorphic association scheme whose underlying set is a nonelementary abelian group and whose strongly regular
graphs have negative Latin square parameters.  To our knowledge, no such example was previously known.
Section~\ref{Hamiltonresult} provides examples of constructions in larger fields, and it also includes a
different proof of a theorem shown in \cite{hamilton}.

\section{A 4-class Amorphic Association Scheme}\label{example}
In this section we use Theorem~\ref{brouwerthm} to construct a 4-class amorphic association scheme in the additive
group of $V=\Ff_4^{2 \ell}$, where $\ell \geq 2$. According to Theorem~\ref{vanDamtheorem}, we can do this by partitioning
$V\setminus\{0\}$ into 4 partial difference sets with negative Latin square type parameters. Let
$\Ff_4=\{0,1,\alpha, \alpha^2=\alpha+1\}$, and let
$$Q(x_1, x_2, \ldots, x_{2 \ell}) = \alpha x_1^2 + x_1 x_2 + x_2^2 +
x_3 x_4 + \cdots + x_{2 \ell - 1} x_{2 \ell}.$$ This is an elliptic quadratic form on $V$. By
Theorem~\ref{ellquadform}, the set $\Omega = \{x\in V\mid Q(x) = 0, x \neq 0\}$ is a negative Latin square type
PDS.  We introduce the notation $$D_{\beta} = \{ v \in \Ff_{4}^{2 \ell}\mid Q(v) = \beta \},$$ where
$\beta\in\Ff_4$. Note that $\Omega = D_0 \backslash \{ 0 \}$ and $D_1=\{v\in \Ff_4^{2\ell}\mid \tr_{4/2} \circ Q
(v)=0\}\setminus D_0$. Theorem~\ref{brouwerthm} shows that the set $D_1$ is a negative Latin square type PDS. Now
we would like to show that $D_{\alpha} = \{x\in V \mid Q(x) = \alpha \}$ and $D_{\alpha^2} = \{x\in V\mid Q(x) =
\alpha^2 \}$ are each negative Latin square type PDSs. We can do this by using the quadratic forms $Q' = \alpha
Q$ and $Q'' = \alpha^2 Q$ and applying Theorems~\ref{ellquadform} and \ref{brouwerthm} to these quadratic forms
in the same way we used $Q$ to show that $D_1$ is a negative Latin square type PDS. Combining these PDSs with
Theorem~\ref{vanDamtheorem}, we have proved the following theorem.

\begin{thm}
\label{specialsituation} The strongly regular Cayley graphs generated by the negative Latin square type PDSs
$D_0\setminus\{0\}$, $D_1$, $D_{\alpha}$, and $D_{\alpha^2}$ form a 4-class amorphic association scheme on $V$.
\end{thm}

\begin{remark}
The subset $D_0\setminus\{0\}$ has size $(4^{\ell}+1)(4^{\ell -1}-1)$, while the subsets $D_1$, $D_{\alpha}$, and
$D_{\alpha^2}$ of $V$ all have the same size, namely $(4^{\ell}+1)4^{\ell -1}$.
\end{remark}

\section{Amorphic association schemes defined on nonelementary abelian 2-groups}
\label{nonelementaryabelianconstruction}

In this section, we demonstrate that we can have an amorphic association scheme on $(G, +)$, in which all graphs
are negative Latin Square type SRGs, and $G$ is a nonelementary abelian 2-group. The construction in this section
follows the techniques found in \cite{davisxiang}, where the set $D_0\setminus\{0\}$ from the previous section
was ``lifted'' from a PDS in $(\Ff_4^{2\ell},+)$ to a PDS in the group $\Z_4^2 \times \Z_2^{4 \ell - 4}$ by using
the map $F$ defined in Section 1.2. We will only include the smallest case where $G=\Z_4^2 \times \Z_2^{4 \ell -
4}$ since the general case (i.e., $G=\Z_4^{2k} \times \Z_2^{4 \ell - 4k}$, $1<k<\ell$) can be handled in a way
completely analogous to that in \cite{davisxiang}.

Again let $\ell\geq 2$ be an integer, $\Ff_4=\{0,1,\alpha,\alpha^2=\alpha +1\}$,
$$Q(x_1, x_2, \ldots, x_{2 \ell}) = \alpha x_1^2 + x_1 x_2 + x_2^2 +
x_3 x_4 + \cdots + x_{2 \ell - 1} x_{2 \ell},$$ and let $D_{\beta} = \{ v \in \Ff_{4}^{2 \ell}\mid Q(v) =
\beta\}$, where $\beta\in\Ff_4$. We define $${\mathcal L}_{\beta} = F(D_{\beta})$$ for every $\beta\in \Ff_4$.
For example,
$${\mathcal L}_1=\{(\xi_1+2\xi_2,x_3,\ldots ,x_{2\ell})\in R\times\Ff_4^{2\ell-2}\mid
Q(\pi(\xi_1),\pi(\xi_2),x_3,\ldots ,x_{2\ell})=1\},$$ where $R=\GR(4,2)$. The goal of this section is to
demonstrate that ${\mathcal L}_{\beta}$, $\beta\in\{1,\alpha,\alpha^2\}$, are all negative Latin square type
PDSs. We will provide the proof that ${\mathcal L}_1$ is a PDS in the additive group of $R\times\Ff_4^{2\ell-2}$
in its full details. It can be shown that ${\mathcal L}_{\alpha}$ and ${\mathcal L}_{\alpha^2}$ are also PDS by
the same arguments with $Q$ replaced by $\alpha Q$ and $\alpha^2 Q$ respectively. The careful reader will notice
that an important part of the proofs (i.e., Lemma~\ref{hardcase}) we are going to give is different from that in
\cite{davisxiang}. When dealing with the ``lifting'' of $D_0\setminus\{0\}$, the geometry of quadrics can help us
with the proof (see Lemma 2.7 of \cite{davisxiang}). It seems difficult to find similar geometric arguments which
also work for ${\mathcal L}_1$. We have to rely on algebraic computations.

By Lemma~\ref{charsum-PDS}, in order to demonstrate that ${\mathcal L}_1$ is a PDS we need to show that
$\Psi({\mathcal L}_1)$ are as required for all nonprincipal characters $\Psi$ of $R\times \Ff_4^{2 \ell - 2}$. If
we write $\Psi = \psi_{\beta_1 + 2 \beta_2} \otimes \chi_{w_3, w_4, \ldots, w_{2 \ell}}$, the character sum of
${\mathcal L}_1$ is as follows.
$$\Psi({\mathcal L}_1) = \Sigma_{Q(\pi(\xi_1),
\pi(\xi_2), x_3, \ldots, x_{2 \ell}) = 1} \sqrt{-1}^{\Tr(\beta_1 \xi_1)} (-1)^{\tr(\pi(\beta_2) \pi(\xi_1) +
\pi(\beta_1) \pi(\xi_2) + w_3 x_3 + \cdots + w_{2 \ell} x_{2 \ell})}$$ If $\beta_1 = 0$, that is, $\Psi$ is a
character of order 2, then $\Psi({\mathcal L}_1) = \chi_{\pi(\beta_2), 0, w_3, \ldots, w_{2 \ell}}(D_1)$.  Since
we know that $D_1$ is a PDS by Theorem~\ref{brouwerthm}, we have $\Psi({\mathcal L}_1) = 4^{\ell - 1}$ or
$4^{\ell-1} - 4^{\ell}$ as required (cf. Lemma~\ref{latinPDS}). Thus, the characters of order 2 all have the
correct sum over ${\mathcal L}_1$. So from now on we assume that $\beta_1 \neq 0$. Then $\Psi({\mathcal L}_1)$
can be broken down into three sub-sums according as $\Tr(\beta_1 \xi_1)$ is $0, 2,$ or odd.
\begin{eqnarray*}
\Psi({\mathcal L}_1)&=&\sum_{(\xi_1+2\xi_2,x_3,\ldots,x_{2\ell})\in {\mathcal L}_1, \Tr(\beta_1\xi_1)=0}(-1)^{\tr(\pi(\beta_2\xi_1+\beta_1\xi_2))}(-1)^{\tr(\sum_{i=3}^{2\ell}w_ix_i)}\\
&{}&-\sum_{(\xi_1+2\xi_2,x_3,\ldots,x_{2\ell})\in {\mathcal L}_1, \Tr(\beta_1\xi_1)=2}(-1)^{\tr(\pi(\beta_2\xi_1+\beta_1\xi_2))}(-1)^{\tr(\sum_{i=3}^{2\ell}w_ix_i)}\\
&{}&+\sum_{(\xi_1+2\xi_2,\xi_3,\ldots,\xi_{2\ell})\in {\mathcal L}_1, \Tr(\beta_1\xi_1)={\rm odd}}(\sqrt{-1})^{\Tr(\beta_1\xi_1)}(-1)^{\tr(\pi(\beta_2\xi_1+\beta_1\xi_2))}(-1)^{\tr(\sum_{i=3}^{2\ell}w_ix_i)}\\
\end{eqnarray*}
The third sub-sum can be shown to be zero in exactly the same way as in \cite{davisxiang}: namely, use the fact
that $Q(x_1, x_1+x_2, x_3, \ldots, x_{2 \ell}) = Q(x_1, x_2, \ldots, x_{2 \ell})$ to identify a pair of elements whose corresponding terms in the third sub-sum add to 0. For the first sub-sum, note that $\Tr(\beta_1 \xi_1) = 0$ implies $\xi_1 = 0$. Define
$$O_0 =\{ (x_1, x_2, \ldots, x_{2 \ell}) \in D_1 \mid x_1 = 0 \}.$$ Then the first sub-sum is equal to
$\chi_{\pi(\beta_2), \pi(\beta_1), w_3, \ldots, w_{2 \ell}}(O_0)$. For the second sub-sum, note that $\Tr(\beta_1\xi_1) = 2$ implies $\xi_1 = \beta_1^{-1}$. We have
\begin{eqnarray*}
\mbox {\rm The second sub-sum}&=&\sum_{(x_1,x_2,\ldots ,x_{2\ell})\in D_1, \; x_1=\pi(\beta_1)^{-1}}\chi_{\pi(\beta_2), \pi(\beta_1), w_3, \ldots, w_{2 \ell}}((x_1,x_2,\ldots ,x_{2\ell}))\\
&=&\sum_{(x_1,x_2,\ldots ,x_{2\ell})\in D_1,\; x_1\neq 0}\chi_{\pi(\beta_2), \pi(\beta_1), w_3, \ldots, w_{2 \ell}}((x_1,x_2,\ldots ,x_{2\ell}))\\
&{}&-\sum_{(x_1,x_2,\ldots ,x_{2\ell})\in D_1,\; x_1\neq 0, \;\tr(\pi(\beta_1)x_1)=1}\chi_{\pi(\beta_2), \pi(\beta_1), w_3, \ldots, w_{2 \ell}}((x_1,x_2,\ldots ,x_{2\ell}))\\
&=&\sum_{(x_1,x_2,\ldots ,x_{2\ell})\in D_1,\; x_1\neq 0}\chi_{\pi(\beta_2), \pi(\beta_1), w_3, \ldots, w_{2 \ell}}((x_1,x_2,\ldots ,x_{2\ell}))\\
\end{eqnarray*}
The last equality follows since $$\sum_{(x_1,x_2,\ldots ,x_{2\ell})\in D_1,\; x_1\neq 0,
\;\tr(\pi(\beta_1)x_1)=1}\chi_{\pi(\beta_2), \pi(\beta_1), w_3, \ldots, w_{2 \ell}}((x_1,x_2,\ldots
,x_{2\ell}))=0$$
by the pairing argument used in dealing with the third sub-sum. Hence
\begin{equation}\label{secondsum}
\mbox{\rm The second sub-sum}=\chi_{\pi(\beta_2), \pi(\beta_1), w_3, \ldots, w_{2 \ell}}(D_1 \backslash O_0)
\end{equation}
Thus,
\begin{eqnarray}\label{ringsum}
\Psi({\mathcal L}_1) & = & \chi_{\pi(\beta_2), \pi(\beta_1), w_3, \ldots, w_{2 \ell}}(O_0) - \chi_{\pi(\beta_2),
\pi(\beta_1), w_3, \ldots, w_{2 \ell}}(D_1 \backslash O_0)
\end{eqnarray}
Contrast this with
\begin{eqnarray}\label{fieldsum}
\chi_{\pi(\beta_2), \pi(\beta_1), w_3, \ldots, w_{2 \ell}}(D_1) & = & \chi_{\pi(\beta_2), \pi(\beta_1), w_3,
\ldots, w_{2 \ell}}(O_0) + \chi_{\pi(\beta_2), \pi(\beta_1), w_3, \ldots, w_{2 \ell}}(D_1 \backslash O_0)
\end{eqnarray}
Plugging $x_1=0$ into $Q(x_1,x_2,\ldots ,x_{2\ell})=1$ yields $x_2^2 + x_3 x_4 + \cdots + x_{2 \ell - 1} x_{2
\ell}=1$, leading to the following simplification:
\begin{eqnarray*}
\chi_{\pi(\beta_2), \pi(\beta_1), w_3, \ldots, w_{2 \ell}}(O_0) &
= & (-1)^{\tr(\pi(\beta_1))} \Sigma_{(x_3, \ldots, x_{2 \ell}) \in
\Ff_4^{2 \ell - 2}} (-1)^{\tr(\pi(\beta_1)^2 (x_3 x_4 + \cdots +
x_{2 \ell -1} x_{2 \ell}) + w_3^2 x_3^2 + \cdots + w_{2 \ell}^2
x_{2 \ell}^2)}
\end{eqnarray*}
Since $\pi(\beta_1)^2 (x_3 x_4 + \cdots + x_{2 \ell -1} x_{2 \ell}) + w_3^2 x_3^2 + \cdots + w_{2 \ell}^2 x_{2
\ell}^2$ is a nonsingular quadratic form, Theorem 3.2 of \cite{leepsch} implies that $$\chi_{\pi(\beta_2),
\pi(\beta_1), w_3, \ldots, w_{2 \ell}}(O_0)=\pm 4^{\ell - 1}$$ In the following two lemmas, we will find the
values of $\chi_{\pi(\beta_2), \pi(\beta_1), w_3, \ldots, w_{2 \ell}}(D_1\setminus O_0)$ according as
$\chi_{\pi(\beta_2), \pi(\beta_1), w_3, \ldots, w_{2 \ell}}(O_0)=-4^{\ell -1}$ or $4^{\ell -1}$.

\begin{lemma}
\label{easycase} If $\chi_{\pi(\beta_2), \pi(\beta_1), w_3,
\ldots, w_{2 \ell}}(O_0) = -4^{\ell - 1}$, then
$\chi_{\pi(\beta_2), \pi(\beta_1), w_3, \ldots, w_{2 \ell}}(D_1
\backslash O_0) = \pm 2 \cdot 4^{\ell - 1}$.
\end{lemma}

\bproof Since $D_1$ is a PDS, the character sum of $D_1$ must be either $4^{\ell-1}$ or $4^{\ell-1} - 4^{\ell}$
(cf. Lemma~\ref{latinPDS}). Note that
$$\chi_{\pi(\beta_2), \pi(\beta_1), w_3, \ldots, w_{2 \ell}}(D_1
\backslash O_0)=\chi_{\pi(\beta_2), \pi(\beta_1), w_3, \ldots, w_{2 \ell}}(D_1)-\chi_{\pi(\beta_2), \pi(\beta_1),
w_3, \ldots, w_{2 \ell}}(O_0).$$ Therefore, if $\chi_{\pi(\beta_2), \pi(\beta_1), w_3, \ldots, w_{2 \ell}}(O_0) =
-4^{\ell - 1}$, then $\chi_{\pi(\beta_2), \pi(\beta_1), w_3, \ldots, w_{2 \ell}}(D_1 \backslash O_0) = \pm 2
\cdot 4^{\ell - 1}$. \eproof

Together with (\ref{ringsum}), the above lemma immediately leads to the following corollary.

\begin{cor}\label{easycase continued}
If $\chi_{\pi(\beta_2), \pi(\beta_1),w_3, \ldots, w_{2 \ell}}(O_0) = -4^{\ell - 1}$, then $\Psi({\mathcal L}_1) =
4^{\ell-1}$ or $4^{\ell-1} - 4^{\ell}$.
\end{cor}

Next we consider the case where $\chi_{\pi(\beta_2), \pi(\beta_1), w_3, \ldots, w_{2 \ell}}(O_0) = 4^{\ell - 1}$.
Again using~(\ref{fieldsum}), we see that in this case $\chi_{\pi(\beta_2), \pi(\beta_1), w_3, \ldots, w_{2
\ell}}(D_1 \backslash O_0) = -4^{\ell}$ or $0$. We will show that the case $\chi_{\pi(\beta_2), \pi(\beta_1),
w_3, \ldots, w_{2 \ell}}(D_1 \backslash O_0) = -4^{\ell}$ can not occur. The proof is somewhat lengthy.

\begin{lemma}
\label{hardcase} If $\chi_{\pi(\beta_2), \pi(\beta_1), w_3, \ldots, w_{2 \ell}}(O_0) = 4^{\ell - 1}$, then
$\chi_{\pi(\beta_2), \pi(\beta_1), w_3, \ldots, w_{2 \ell}}(D_1 \backslash O_0) = 0$.
\end{lemma}

\bproof From (\ref{secondsum}), we see that the character sum over $D_1 \backslash O_0$ is equal to the character
sum over the elements of $D_1$ that satisfy $x_1 = \pi(\beta_1)^{-1}$. Plugging this condition into the equation
$Q(x_1,x_2,\ldots ,x_{2\ell})=1$ yields
\begin{equation}\label{defequ}
\pi(\beta_1)x_2+\pi(\beta_1)^2x_2^2+\alpha+\pi(\beta_1)^2+\pi(\beta_1)^2(x_3x_4+\cdots +x_{2\ell-1}x_{2\ell})=0
\end{equation}
For convenience, set
$$\gamma:= \alpha + \pi(\beta_1)^2 + \pi(\beta_1)^2 (x_3 x_4 + \cdots +
x_{2 \ell - 1} x_{2 \ell})$$ Using (\ref{defequ}) in the sum $\chi_{\pi(\beta_2), \pi(\beta_1), w_3, \ldots,
w_{2\ell}}(D_1 \backslash O_0)$, and using properties of the finite field trace, we have
\begin{equation}\label{sumDO}
\chi_{\pi(\beta_2), \pi(\beta_1), w_3, \ldots, w_{2 \ell}}(D_1 \backslash O_0) = (-1)^{\tr(\pi(\beta_2)
\pi(\beta_1)^{-1})} \cdot 2 \cdot (S_1 - S_2),
\end{equation}
where
$$S_1=\Sigma_{(x_3, x_4, \ldots, x_{2 \ell}) \in \Ff_4^{2 \ell - 2},\;
\gamma=0} (-1)^{\tr(w_3 x_3 + \cdots + w_{2 \ell} x_{2\ell})},$$
and
$$S_2=\Sigma_{(x_3, x_4, \ldots, x_{2 \ell}) \in \Ff_4^{2 \ell - 2},\; \gamma=1} (-1)^{\tr(w_3 x_3 + \cdots
+ w_{2 \ell} x_{2 \ell})}.$$  Following a similar computation over the elements $(0,x_2,\ldots ,x_{2\ell})\in
O_0$ and defining $\gamma$ as above (here $x_2^2=1+x_3x_4+\cdots +x_{2\ell -1}x_{2\ell}$ implies that $\gamma=\alpha + \pi(\beta_1)^2x_2^2$), we see that
\begin{eqnarray}\label{sumO_0}
\chi_{\pi(\beta_2), \pi(\beta_1), w_3, \ldots, w_{2 \ell}}(O_0) & = & \Sigma_{(x_3, x_4, \ldots, x_{2 \ell}) \in
\Ff_4^{2 \ell - 2},\; \gamma\in \{ \alpha, \alpha^2 \}} (-1)^{\tr(w_3 x_3 + \cdots + w_{2 \ell} x_{2 \ell})}\nonumber \\
& - & \Sigma_{(x_3, x_4, \ldots, x_{2 \ell}) \in \Ff_4^{2 \ell - 2},\; \gamma\in \{ 0,1 \}} (-1)^{\tr(w_3 x_3 +
\cdots + w_{2 \ell} x_{2 \ell})}
\end{eqnarray}
The main difference between the sum over $D_1\setminus O_0$ and the sum over $O_0$ is that in the former case we
know that $\gamma$ is either 0 or 1 because it is equal to $\pi(\beta_1)^2 x_2^2 + \pi(\beta_1)
x_2=\tr(\pi(\beta_1) x_2)$ by (\ref{defequ}); while in the latter case, $\gamma$ can be any element in $\Ff_4$.
For convenience, we will call the first sum in (\ref{sumO_0}) $S_3$, and the second sum $S_4$. It follows from
the definition of $S_1$, $S_2$, and $S_4$ that
$$S_1+S_2=S_4$$
By assumption, $\chi_{\pi(\beta_2), \pi(\beta_1), w_3, \ldots, w_{2 \ell}}(O_0) = 4^{\ell-1}$, we have
$$S_3-S_4=4^{\ell-1}$$
If we add $S_3$ and $S_4$, we are simply taking the character sum over all of $\Ff_4^{2 \ell - 2}$. So
\begin{equation}
S_3+S_4= \left\{
\begin{array}{ll}
4^{2\ell -2}, & \mbox{if} \; (w_3,w_4,\ldots ,w_{2\ell})=(0,0,\ldots ,0),\\
0, & \mbox{otherwise}. \end{array} \right.
\end{equation}
Hence we have
\begin{equation}\label{assumpS_1+S_2}
2(S_1+S_2)=2S_4= \left\{
\begin{array}{ll}
4^{2\ell -2}-4^{\ell -1}, & \mbox{if} \; (w_3,w_4,\ldots ,w_{2\ell})=(0,0,\ldots ,0),\\
-4^{\ell -1}, & \mbox{otherwise}. \end{array} \right.
\end{equation}

We now return to the computation of $\chi_{\pi(\beta_2), \pi(\beta_1), w_3, \ldots, w_{2 \ell}}(D_1 \backslash
O_0)$. Define
$$H(x_3,x_4,\ldots ,x_{2\ell})=\pi(\beta_1)^2(x_3x_4+\cdots +x_{2\ell -1}x_{2\ell}).$$
Note that $H$ is a nonsingular hyperbolic quadratic form on $\Ff_4^{2\ell -2}$, and
\begin{eqnarray*}
S_1 &=& \sum_{H(x_3,\ldots ,x_{2\ell})=\alpha+\pi(\beta_1)^2}\chi_{w_3,\ldots ,w_{2\ell}}((x_3,x_4,\ldots
,x_{2\ell})),\\
S_2 &=& \sum_{H(x_3,\ldots ,x_{2\ell})=1+\alpha+\pi(\beta_1)^2}\chi_{w_3,\ldots ,w_{2\ell}}((x_3,x_4,\ldots
,x_{2\ell})).
\end{eqnarray*}
The sums $S_1$ and $S_2$ can be evaluated. If $\alpha +\pi(\beta_1)^2=0$, then in the sum for $S_1$, we are
summing over all $(x_3,\ldots ,x_{2\ell})$ satisfying $H(x_3,\ldots ,x_{2\ell})=0$, including the element
$(0,0,\ldots ,0)$. So
$$S_1=1+\sum_{H(x_3,\ldots ,x_{2\ell})=0, (x_3,\ldots ,x_{2\ell})\neq (0,0,\ldots 0)}\chi_{w_3,\ldots ,w_{2\ell}}((x_3,x_4,\ldots
,x_{2\ell})).$$ Now note that the subset $\{(x_3,\ldots ,x_{2\ell})\mid H(x_3,\ldots ,x_{2\ell})=0, (x_3,\ldots
,x_{2\ell})\neq (0,0,\ldots 0)\}$ is a Latin square type PDS in $\Ff_4^{2\ell -2}$ (this follows from the
hyperbolic version of Theorem~\ref{ellquadform}). By Lemma~\ref{latinPDS}, we have
\begin{equation}\label{s1case1}
S_1= \left\{
\begin{array}{ll}
4^{2\ell -3}+4^{\ell -1}-4^{\ell -2}, & \mbox{if} \; (w_3,w_4,\ldots ,w_{2\ell})=(0,0,\ldots ,0),\\
4^{\ell -1}-4^{\ell -2} \; {\mbox or}\; -4^{\ell -2}, & \mbox{otherwise}. \end{array} \right.
\end{equation}
If $\alpha +\pi(\beta_1)^2\neq 0$, then by Theorem~\ref{brouwerthm} (with the quadratic form $Q$ in that theorem
being $H$), we have
\begin{equation}\label{s1case2}
S_1= \left\{
\begin{array}{ll}
4^{2\ell -3}-4^{\ell -2}, & \mbox{if} \; (w_3,w_4,\ldots ,w_{2\ell})=(0,0,\ldots ,0),\\
4^{\ell -1}-4^{\ell -2} \; {\mbox or}\; -4^{\ell -2}, & \mbox{otherwise}. \end{array} \right.
\end{equation}
The same is true for $S_2$. That is, if $1+\alpha +\pi(\beta_1)^2=0$, then
\begin{equation}\label{s2case1}
S_2= \left\{
\begin{array}{ll}
4^{2\ell -3}+4^{\ell -1}-4^{\ell -2}, & \mbox{if} \; (w_3,w_4,\ldots ,w_{2\ell})=(0,0,\ldots ,0),\\
4^{\ell -1}-4^{\ell -2} \; {\mbox or}\; -4^{\ell -2}, & \mbox{otherwise}. \end{array} \right.
\end{equation}
And if $1+\alpha +\pi(\beta_1)^2\neq 0$,
\begin{equation}\label{s2case2}
S_2= \left\{
\begin{array}{ll}
4^{2\ell -3}-4^{\ell -2}, & \mbox{if} \; (w_3,w_4,\ldots ,w_{2\ell})=(0,0,\ldots ,0),\\
4^{\ell -1}-4^{\ell -2} \; {\mbox or}\; -4^{\ell -2}, & \mbox{otherwise}. \end{array} \right.
\end{equation}

If $(w_3,w_4,\ldots ,w_{2\ell})=(0,0,\ldots ,0)$, we contend that $S_1=S_2=4^{2\ell -3}-4^{\ell -2}$. The reason
is that all other choices for $S_1$ and $S_2$ result in
$$S_1+S_2=(4^{2\ell -3}+4^{\ell -1}-4^{\ell -2})+(4^{2\ell -3}-4^{\ell -2})=2(4^{2\ell -3}-4^{\ell -2})+4^{\ell -1},$$
which is never equal to $(4^{2\ell -2}-4^{\ell -1})/2$ as required by (\ref{assumpS_1+S_2}). Hence in this case
$\chi_{\pi(\beta_2), \pi(\beta_1), w_3, \ldots, w_{2 \ell}}(D_1 \backslash O_0)=0$ by (\ref{sumDO}).

If $(w_3,w_4,\ldots ,w_{2\ell})\neq (0,0,\ldots ,0)$, we contend that $S_1=S_2=-4^{\ell-2}$. The reason is that
all other choices for $S_1$ and $S_2$ result in either
$$S_1+S_2=4^{\ell -1} - 4^{\ell -2}+(-4^{\ell -2})=2\cdot 4^{\ell -2},$$
or
$$S_1+S_2=2(4^{\ell -1}-4^{\ell -2})=6\cdot 4^{\ell -2}.$$
Neither of these two values equals $-2\cdot 4^{\ell -2}$ as required by (\ref{assumpS_1+S_2}). Hence in this case
we also have $\chi_{\pi(\beta_2), \pi(\beta_1), w_3, \ldots, w_{2 \ell}}(D_1 \backslash O_0)=0$ by (\ref{sumDO}).
This completes the proof of the lemma. \eproof

\begin{cor}\label{hardcase continued}
If $\chi_{\pi(\beta_2), \pi(\beta_1),w_3, \ldots, w_{2 \ell}}(O_0) = 4^{\ell - 1}$, then $\Psi({\mathcal L}_1) = 4^{\ell-1}$.
\end{cor}

\bproof Use (\ref{ringsum}) and Lemma~\ref{hardcase}.\eproof

Combining the discussion on characters of order 2, Corollary~\ref{easycase continued} and Corollary~\ref{hardcase
continued}, we have the following theorem.

\begin{thm}
\label{liftisaPDS} The set ${\mathcal L}_1$ is a negative Latin
square type PDS in the additive group of $R\times \Ff_4^{2
\ell - 2}$.
\end{thm}

\bproof Let $\Psi$ be a character of the additive group of $R\times \Ff_4^{2 \ell - 2}$. As we see in the above,
the character sum $\Psi({\mathcal L}_1)$ has the correct values in all cases. The theorem follows from
Lemma~\ref{latinPDS}. \eproof

As in Section~\ref{example}, the arguments to show that ${\mathcal L}_{\alpha}$ and ${\mathcal L}_{\alpha^2}$ are
negative Latin square type PDSs will go through in analogous ways by considering the quadratic forms $\alpha Q$
and $\alpha^2 Q$.  The subset ${\mathcal L}_0\setminus\{0\}$ was already shown in \cite{davisxiang} to be a
negative Latin square type PDS. This leads to the main result of this section.

\begin{thm}
\label{maintheoremnonelementaryabelian} The strongly regular Cayley graphs associated to the negative Latin
square type PDSs ${\mathcal L}_0\setminus\{0\}, {\mathcal L}_1, {\mathcal L}_{\alpha},$ and
${\mathcal L}_{\alpha^2}$ form a 4-class amorphic association scheme on $R\times \Ff_4^{2\ell -2}$.
\end{thm}

\bproof  Note that ${\mathcal L}_0\setminus\{0\}, {\mathcal L}_1, {\mathcal L}_{\alpha},$ and ${\mathcal L}_{\alpha^2}$ partition $R\times \Ff_4^{2\ell -2}\setminus\{0\}$. Each of them is a negative Latin square type PDS. By Theorem~\ref{vanDamtheorem}, the strongly regular Cayley graphs associated with these four subsets form an amorphic association scheme on $R\times \Ff_4^{2\ell -2}$. \eproof

As indicated in the introduction, this association scheme seems to be the first amorphic association scheme with
negative Latin square type graphs defined on a nonelementary abelian group. (We remark that there are known
constructions of amorphic association schemes with Latin square type graphs defined on the additive groups of
Galois rings. See \cite{imy}.) We could extend this result to more copies of $R$ using techniques similar to that
in \cite{davisxiang}, but the computations are tedious. We therefore will not include the details here. The
interested reader is referred to \cite{davisxiang} for details on how to construct PDSs in $R^{k}\times
\Ff_4^{2\ell -2k}$, where $1<k<\ell$.

\section{Further amorphic schemes from Brouwer's construction}
\label{Hamiltonresult}

Let $q$ be a prime power, and let $m, \ell$ be positive integers. Suppose that $m$ has a chain of divisors
$$1=m_d|m_{d-1}|m_{d-2}|\cdots |m_2|m_1=m,$$
where $m_{i}\neq m_{i-1}$ for all $1<i\leq d$. Then the finite field $\Ff_{q^m}$ has a chain of subfields
$$\Ff_q=\Ff_{q^{m_d}}\subset \Ff_{q^{m_{d-1}}} \subset \cdots \subset \Ff_{q^{m_2}}\subset \Ff_{q^{m_1}} = \Ff_{q^m}.$$
Let $Q: V\rightarrow \Ff_{q^m}$ be a nonsingular elliptic quadratic form, where $V=\Ff_{q^m}^{2\ell}$, and let
$\tr_{q^m/q^{m_i}}$ be the trace from $\Ff_{q^m}$ to $\Ff_{q^{m_i}}$. For $1\leq i\leq d$, we define
$$Q_i=\tr_{q^m/q^{m_i}}\circ Q,$$
and
$$\Omega_i=\{x\in V\mid x\neq 0, Q_i(x)=0\}.$$
Then each $Q_i: V\rightarrow \Ff_{q^{m_i}}$ is a nonsingular elliptic quadratic form (cf. \cite{brouwer}), and by
transitivity of the traces, we have
$$\Omega_1\subset \Omega_2\subset \cdots \subset \Omega_{d}.$$
By Theorem~\ref{ellquadform} each $\Omega_i$ is a negative Latin square type PDS, and Theorem~\ref{brouwerthm}
shows that $\Omega_i \backslash \Omega_j$ is a negative Latin square type PDS for all $i
> j$.  We can prove the following theorem:

\begin{thm}
\label{hamilton'sassociation scheme} Let $Q: V\rightarrow \Ff_{q^m}$ be a nonsingular elliptic quadratic form,
where $V=\Ff_{q^m}^{2\ell}$, and let $\Omega_1 \subset \Omega_2 \subset \cdots \subset \Omega_{d}$ be the sets
described above. Then the Cayley graphs generated by the collection $\Omega_1, \Omega_2 \backslash \Omega_1,
\Omega_3 \backslash \Omega_2, \ldots, \Omega_{d} \backslash \Omega_{d-1}, (V\setminus\{(0,0,\ldots
,0)\})\backslash \Omega_{d}$ form a $(d+1)$-class amorphic association scheme on $V$.
\end{thm}

\bproof Theorem~\ref{brouwerthm} implies that $\Omega_i \backslash \Omega_{i-1}$, $1 < i \leq d$, are negative
Latin square type PDSs when $Q$ is elliptic.  In addition, the complement of $\Omega_d$ in
$V\setminus\{(0,0,\ldots ,0)\}$ is also a PDS of the same type. Thus, we have partitioned
$V\setminus\{(0,0,\ldots ,0)\}$ into disjoint PDSs with negative Latin square type parameters. The strongly
regular Cayley graphs associated with $\Omega_1, \Omega_2 \backslash \Omega_1, \Omega_3 \backslash \Omega_2,
\ldots, \Omega_{d} \backslash \Omega_{d-1}, (V\setminus\{(0,0,\ldots ,0)\})\backslash \Omega_{d}$ are all of
negative Latin square type, and form an edge disjoint decomposition of the complete graph on the vertex set $V$.
Theorem~\ref{vanDamtheorem} implies that these must be the classes of an amorphic association scheme. \eproof

This theorem also works for hyperbolic quadratic forms, producing Latin square type strongly regular graphs as
the classes of the amorphic association scheme.  We note that this theorem contains some ambiguity regarding the
choice of the subfield chain.  For example, if we are working with the field $\Ff_{p^{30}}$, where $p$ is a prime, we
could use any of the following chains: $\Ff_p \subset \Ff_{p^3} \subset \Ff_{p^6} \subset \Ff_{p^{30}}; \Ff_p
\subset \Ff_{p^3} \subset \Ff_{p^{15}} \subset \Ff_{p^{30}}; \Ff_p \subset \Ff_{p^5} \subset \Ff_{p^{10}} \subset
\Ff_{p^{30}}; \Ff_p \subset \Ff_{p^5} \subset \Ff_{p^{15}} \subset \Ff_{p^{30}}; \Ff_p \subset \Ff_{p^2} \subset
\Ff_{p^6} \subset \Ff_{p^{30}}; \Ff_p \subset \Ff_{p^2} \subset \Ff_{p^{10}} \subset \Ff_{p^{30}}.$
Theorem~\ref{hamilton'sassociation scheme} applies to each of these chains to produce a 5-class amorphic
association scheme, but the actual strongly regular graphs will be different for each choice.

Theorem~\ref{hamilton'sassociation scheme} provides an alternative proof of the following result originally found
in Hamilton's paper~\cite{hamilton}.

\begin{cor}
\label{hamilton'sresult} Let $Q: V\rightarrow \Ff_{q^m}$ be a nonsingular elliptic quadratic form, where
$V=\Ff_{q^m}^{2\ell}$, and let $\Omega_1 \subset \Omega_2 \subset \cdots \subset \Omega_{d}$ be described as
above. Then $(\Omega_{d} \backslash \Omega_{d-1}) \cup (\Omega_{d-2} \backslash \Omega_{d-3}) \cup \cdots \cup
(\Omega_{2} \backslash \Omega_{1})$ is a negative Latin square type PDS if $d$ is even and $(\Omega_{d}
\backslash \Omega_{d-1}) \cup (\Omega_{d-2} \backslash \Omega_{d-3}) \cup \cdots \cup (\Omega_{3} \backslash
\Omega_{2}) \cup \Omega_{1}$ is a negative Latin square type PDS if $d$ is odd.
\end{cor}

\bproof  Simply fuse the appropriate classes from the amorphic
association scheme in Theorem~\ref{hamilton'sassociation scheme}.
\eproof

The simple proof of Corollary~\ref{hamilton'sresult} given above demonstrates the power of
Theorem~\ref{vanDamtheorem}. We show one more application of this theorem by doing a variation of
Theorem~\ref{hamilton'sassociation scheme} for quadratic forms with values in $\Ff_{q^m}$, where $m$ is even.
Without loss of generality, we will assume $m=2$. The result will be a generalization of
Theorem~\ref{specialsituation} in Section~\ref{example}. Before we state the result, we need the following lemma
identifying the trace 0 elements in a quadratic extension field of $\Ff_q$.

\begin{lemma}
\label{quadraticextensiontracezeroelements}
\begin{enumerate}
\item Let $s\geq 1$ be an integer, $q = 2^s$, and let $g$ be a primitive element of $\Ff_{q^2}$.  Then $\{ x \in \Ff_{q^2} \mid \tr_{q^2/q}(x) = 0 \} = \{ 0 \} \cup \{ g^{(q+1)i} \mid 0 \leq i <
q-1 \}$.
\item Let $s\geq 1$ be an integer, $p$ an odd prime, $q = p^s$, and let $g$ be a primitive element of $\Ff_{q^2}$.  Then $\{ x \in \Ff_{q^2} \mid \tr_{q^2/q}(x) = 0 \} = \{ 0 \}
\cup \{ g^{\frac{q+1}{2} + (q+1)i} \mid 0 \leq i < q-1 \}$.
\end{enumerate}
\end{lemma}

\bproof In both cases, we have $\tr_{q^2/q}(x) = x + x^q$ for all $x\in \Ff_{q^2}$. In the $q=2^s$ case, note
that $g^{(q+1)i}\in \Ff_q\subset \Ff_{q^2}$, hence $\tr_{q^2/q}(g^{(q+1)i}) = g^{(q+1)i} + (g^{(q+1)i})^q =
g^{(q+1)i} + g^{(q+1)i} = 0$. In the $q = p^s$ case, where $p$ is an odd prime, again noting that $g^{(q+1)i}\in
\Ff_q\subset \Ff_{q^2}$, we have
$$\tr_{q^2/q}(g^{\frac{q+1}{2} + (q+1)i}) = g^{(q+1)i}(g^{\frac{q+1}{2}} + g^{\frac{(q+1)q}{2}})=g^{(q+1)i}(g^{\frac{q+1}{2}}+g^{\frac{q^2-1}{2}}g^{\frac{q+1}{2}})=0,$$
since $g^{\frac{q^2-1}{2}}=-1$. A counting argument indicates that we have found all of the trace 0 elements.
\eproof

We now vary the quadratic form in a way similar to our approach in Section~\ref{example}: if
$Q:\Ff_{q^2}^{2\ell}\rightarrow \Ff_{q^2}$ is a nonsingular elliptic quadratic form and $g$ is a primitive
element of $\Ff_{q^2}$, then consider the $q+1$ quadratic forms $\tr_{q^2/q} \circ g^iQ, 0 \leq i \leq q$. Let
$$\Omega_{g^i} = \{ x \in \Ff_{q^2}^{2 \ell} \mid \tr_{q^2/q} \circ g^iQ(x) = 0, x \neq 0 \},$$
for all $0\leq i\leq q$. Each of these sets contains the same subset $\Omega_0 = \{ x \in \Ff_{q^2}^{2 \ell} \mid
x\neq 0, Q(x) = 0 \}$ associated to the quadratic form $Q$ since $Q(x) = 0$ if and only if $g^iQ(x) = 0$.
Theorem~\ref{brouwerthm} implies that $\Omega_{g^i} \backslash \Omega_0$ is a negative Latin square type PDS for
all $i=0,1,\ldots ,q$. This leads to the following theorem.

\begin{thm}
\label{cyclicrotationofquadrics} Let $g$ be a primitive element of $\Ff_{q^2}$, and let
$Q:\Ff_{q^2}^{2\ell}\rightarrow \Ff_{q^2}$ be a nonsingular elliptic quadratic form.  Then the strongly regular
graphs associated to the PDSs $\Omega_0, \Omega_{g^i}\backslash \Omega_0, 0 \leq i \leq q$, form a $(q+2)$-class
amorphic association scheme.
\end{thm}

\bproof Each of the PDSs $\Omega_0, \Omega_{g^i}\backslash \Omega_0, 0 \leq i \leq q$, is of negative Latin
square type, and hence the associated strongly regular Cayley graph is also of negative Latin square type.
Lemma~\ref{quadraticextensiontracezeroelements} implies that these negative Latin square type PDSs partition
$\Ff_{q^2}^{2 \ell}\setminus\{0\}$. That means the associated SRGs form an edge-decomposition of the complete
graph on $\Ff_{q^2}^{2\ell}$. Now the result follows from Theorem~\ref{vanDamtheorem}. \eproof

As an example, consider a nonsingular elliptic quadratic form with its values in $\Ff_{p^{30}}$. If we apply
Theorem~\ref{hamilton'sassociation scheme}, we will get several 5-class amorphic association schemes. If we apply
Theorem~\ref{cyclicrotationofquadrics} with $q=p^{15}$, we get a $(p^{15}+2)$-class amorphic association scheme.
These can be fused to yield any of the 5-class amorphic association schemes from the previous construction.

Setting $q=2$ in Theorem~\ref{cyclicrotationofquadrics}, we obtain Theorem~\ref{specialsituation}. In
Section~\ref{nonelementaryabelianconstruction}, we ``lifted'' the PDSs $D_0\setminus\{0\}$, $D_1$, $D_{\alpha}$,
and $D_{\alpha^2}$ in Theorem~\ref{specialsituation} to PDSs ${\mathcal L}_0\setminus\{0\}$, ${\mathcal L}_1$,
${\mathcal L}_{\alpha}$, and ${\mathcal L}_{\alpha^2}$ in $(R\times \Ff_4^{2\ell -2}, +)$ that partition
$(R\times \Ff_4^{2\ell -2})\setminus\{(0,0,\ldots ,0)\}$, hence obtained a 4-class amorphic association scheme
defined on $(R\times \Ff_4^{2\ell -2}, +)$. It is natural to ask whether one can do the same for the PDSs in
Theorem~\ref{cyclicrotationofquadrics} when $q>2$. Computer computations indicate that when $q>2$ and $q=p^s$,
$p$ a prime, if we follow the same procedure to ``lift'' $\Omega_0$ to $R'\times \Ff_{q^2}^{2\ell -2}$, where
$R'$ is the Galois ring $\GR(p^2, 2s)$, the resulting set is not a PDS.

\vspace{0.1in}

\noindent{\bf Acknowledgment:} We thank Frank Fiedler for some initial computer computations. The second author
thanks Department of Mathematics and Computer Science, University of Richmond for the opportunity to spend part
of his sabbatical in March 2004. The first author was supported in part by NSA grant MDA904-03-1-0032. The second
author was supported in part by NSA grant MDA904-03-1-0095.


\begin{thebibliography}{10}

\bibitem{bmw} L. D. Baumert, W. H. Mills, and R. L. Ward, Uniform cyclotomy, {\em J. Number Theory} {\bf 14} (1982), 67--82.

\bibitem{bcn} A. E. Brouwer, A. M. Cohen, and A. Neumaier, {\em Distance Regular Graphs}, Ergebnisse der Mathematik und ihrer Grenzgebiete (3) [Results in Mathematics and Related Areas (3)], 18. Springer-Verlag, Berlin, 1989.

\bibitem{brouwer} A. E. Brouwer, Some new two-weight codes and
strongly regular graphs, {\it Discrete Appl. Math.} {\bf 10}
(1985), 111--114.

\bibitem{calderbankkantor} R. A. Calderbank, W. M. Kantor, The geometry of two-weight codes, {\it Bull. London
Math. Soc.} {\bf 18} (1986), 97--122.

\bibitem{cameron} P. J. Cameron, Finite geometry and coding theory, Lecture Notes for Socrates Intensive Programme ``Finite Geometries and Their Automorphisms'', Potenza, Italy, June 1999. See http://dwispc8.vub.ac.be/Potenza/cameron.ps

\bibitem{vanDam} E. R. van Dam, Strongly regular decompositions of the complete graph, {\em J. Alg. Comb.} {\bf 17} (2003), 181--201.

\bibitem{davisxiang} J. A. Davis, Q. Xiang, Negative Latin square type partial difference sets in nonelementary abelian 2-groups, {\em Journal London Math. Soc.} (2) {\bf 70} (2004), 125--141.

\bibitem{dels} P. Delsarte, An algebraic approach to the association schemes of coding theory, {\em Philips Research Report}, Suppl. No. 10.

\bibitem{hamilton} N. Hamilton, Strongly regular graphs from differences of quadrics, {\em Discrete Math.} {\bf 256} (2002),
465--469.


\bibitem{imy} T. Ito, A. Munemasa, and M. Yamada, Amorphous association schemes over the Galois rings of characteristic 4,
{\it Europ. J. Combin.} {\bf 12} (1991), 513--526.

\bibitem{ivanov} A. V. Ivanov, Amorphous celluar rings II, in: {\it Investigations in Algebraic Theory of Combinatorial Objects}, VNIISI, Moscow, Institute for System Studies, 1985, 39--49 (in Russian).


\bibitem{leepsch} D. Leep,  L. M. Schueller, Zeros of a pair of quadratic forms defined over a finite field,  {\it Finite Fields Appl.}  {\bf 5}  (1999), 157--176.


\bibitem{ln} R. Lidl and H. Niederreiter, {\em Finite Fields}, Cambridge University Press, Cambridge, 1997.


\bibitem{vanLintWilson} J. H. van Lint, R. M. Wilson, {\em A Course in Combinatorics}, Second edition, Cambridge University Press, 2001.


\bibitem{masurvey} S.~L. Ma, A survey of partial difference sets, {\it Designs,
Codes and Cryptography} {\bf 4} (1994), 221--261.

\bibitem{mcdonald} B. R. McDonald, {\it Finite rings with identity}, Pure and Applied Mathematics, Vol. 28. Marcel Dekker, Inc., New York, 1974.

\bibitem{turyn} R.~J. Turyn,\newblock Character sums and difference sets, \newblock {\em Pacific J. Math.} {\bf 15} (1965), 319--346.




\end{thebibliography}
\end{document}